\documentclass[12pt]{amsart}
\usepackage{amsmath,amssymb, euscript, graphicx}

\theoremstyle{definition}

\begin{document}

\address{Azer Akhmedov, Department of Mathematics,
North Dakota State University,
Fargo, ND, 58108, USA}
\email{azer.akhmedov@ndsu.edu}

\begin{center} {\bf EXTENSION OF H\"OLDER'S THEOREM IN Diff $_{+}^{1+\epsilon }(I)$} \end{center}

\medskip

 \begin{center} Azer Akhmedov \end{center}

\medskip

{\bf Abstract:} {\Small We prove that if $\Gamma $ is subgroup of Diff $_{+}^{1+\epsilon }(I)$ and $N$ is a natural number such that every non-identity element of $\Gamma $ has at most $N$ fixed points then $\Gamma $ is solvable. If in addition $\Gamma $ is a subgroup of Diff $_{+}^{2}(I)$ then we can claim that $\Gamma $ is metaabelian.} 

\vspace{1cm}

  It is a classical result (essentially due to H\"older,  cf.[N1]) that if $\Gamma $ is a subgroup of Homeo$_{+}(\mathbb{R})$ such that every nontrivial element acts freely then $\Gamma $ is Abelian. A natural question to ask is what if every nontrivial element has at most $N$ fixed points where $N$ is a fixed natural number. In the case of $N=1$, we do have a complete answer to this question: it has been proved by Solodov (not published), Barbot [B], and Kovacevic [K] that in this case the group is metaabelian, in fact, it is isomorphic to a subgroup of the affine group Aff$(\mathbb{R})$. (see [FF] for the history of this result, where yet another nice proof is presented). 
  
 \medskip
 
 In this paper, we answer this question for an arbitrary $N$ assuming some regularity on the action of the group. 
  
 \medskip
 
  Our main result is the following theorem. 
  
  \medskip
  
  {\bf Theorem 0.1.(Main Theorem)} Let $\epsilon \in (0,1)$ and $\Gamma $ be a subgroup of Diff $_{+}^{1+\epsilon }(I)$ such that every nontrivial element of $\Gamma $ has at most $N$ fixed points. Then $\Gamma $ is solvable.
  
  \medskip
  
  Assuming a higher regularity on the action we obtain a stronger result
  
  \medskip
   
  {\bf Theorem 0.2.} Let $\Gamma $ be a subgroup of Diff $_{+}^{2}(I)$ such that every nontrivial element of $\Gamma $ has at most $N$ fixed points. Then $\Gamma $ is metaabelian.
   
  \medskip
  
  An important tool in obtaining these results is provided by Theorems B-C from [A]. Theorem B (Theorem C) states that a non-solvable (non-metaabelian) subgroup of Diff $_{+}^{1+\epsilon }(I)$ (of Diff $_{+}^{2}(I)$) is non-discrete in $C_0$ metric. Existence of $C_0$-small elements in a group provides effective tools in tackling the problem. Such tools are absent for less regular actions, and the problem for Homeo$_{+}(I)$ (even for Diff$_{+}(I)$) still remains open.

  \medskip   
  
  {\bf Basic Notations:} Throughout this paper, $G$ will denote the group 
 Diff $_{+}^{1+\epsilon }(I)$ where $\epsilon \in (0,1)$. We let $\Gamma \leq G$. For every $g\in \Gamma $, $Fix(g)$ will denote the set of fixed of points of $g$ in $(0,1)$. A fixed point $x_0\in Fix(g)$ is called {\em tangential} if the function $h(x) = g(x)-x$ does not change its sign near $x_0$. For $f,g\in G$, we will say that (the graphs of) $f$ and $g$ have {\em a crossing} at $x_0\in (0,1)$ if $f(x_0) = g(x_0)$. We will write $F(\Gamma ) = \displaystyle \sup_{g\in \Gamma \backslash \{1\}} |Fix(g)|$. (so $F(\Gamma )$ is either 0, or a positive integer, or infinity.) For $f\in G$, we will write $$d(f) = \displaystyle \sup _{x\in [0,1]}|f(x)-x| \ \mathrm{and} \  s(f) = \displaystyle \inf _{x, y\in Fix(f)\cup \{0,1\}, x\neq y}|x-y|$$  
  
  \medskip
  
  If $F(\Gamma ) < \infty $, then we can introduce the following natural biorder in $\Gamma $: for $f, g\in \Gamma $, we write $g < f$ if $g(x) < f(x)$ near zero. If $f$ is a positive element (w.r.t. the order) then we will also write $g<<f$ if $g^n<f$ for every integer $n$; we will say that $g$ is {\em infinitesimal} w.r.t. $f$. For $f\in \Gamma $, we write $\Gamma _f = \{\gamma \in \Gamma : \gamma << f \}$ (so $\Gamma _{f}$ consists of diffeomorphisms which are infinitesimal w.r.t. $f$). Notice that if $\Gamma $ is finitely generated with a fixed finite symmetric  generating set, and $f$ is the biggest generator, then $\Gamma _f$ is a normal subgroup of $\Gamma $, moreover, $\Gamma /\Gamma _f$ is Archimedean, hence Abelian. Notice also that $[\Gamma , \Gamma ]\leq \Gamma _{f}$.

  \medskip

  For the rest of the paper we may and will assume that the action of $\Gamma $ is irreducible, i.e. $\Gamma $ has no global fixed point in $(0,1)$. All the used generating sets of groups will be assumed to be symmetric. Let us also clarify that we define a metaabelian group as a solvable group of derived length at most two (so, Abelian groups are special cases of metaabelian groups). 
  
  \medskip
  
  {\em Acknowledgment.} We are grateful to Andr\' es Navas for very useful conversations and for his valuable remarks on the manuscript.

  \vspace{1cm}
  
  \begin{center} {\bf 1. Preliminary Results} \end{center}
  
  \bigskip
  
  First, we will discuss the case of $N=1$. Let us first show a quick proof of the fact that $N=1\Rightarrow $ {\em metaabelian}.
    
  \medskip
  
  {\bf Proposition 1.1.} If $F(\Gamma ) = 1$ then $\Gamma $ is metaabelian. 
  
  \medskip
  
  {\bf Proof.} (following Farb-Franks, [FF]) If all finitely generated subgroups of a group are metaabelian then the group is metaabelian. Hence we may assume that $\Gamma $ is finitely generated with a fixed finite generating set. 
  
   Let $f$ be the biggest generator of $\Gamma $. Let also  $h\in \Gamma _{f}\backslash \{1\}$ such that $h$ has at least one fixed point (if such $h$ does not exists then $\Gamma _{f}$ is Abelian, therefore $\Gamma $ is metaabelian) and $h(x) > x$ near zero. We may also assume that $\Gamma _f$ has no global fixed point. (if $\Gamma _f$ has a global fixed point then it is Abelian, hence $\Gamma $ is metaabelian.)
 
 \medskip
    
   If $f$ has no fixed point, then for sufficiently big $n$, $fh^{-n}$ will have at least two fixed points, hence, contradiction.
   
   \medskip
   
 Assume now $f$ has one fixed point. Let $a$ be the fixed point of $h$. By conjugating $f$ by the element of $\Gamma _f$ if necessary, we may assume that the fixed point of $f$ is bigger than $a$. Then, again, then for sufficiently big $n$, $fh^{-n}$ will have at least two fixed points.   $\square $  
   
   \medskip
   
   Notice that the above proof works also in Homeo $_{+}(I)$.
   
   \medskip
  
  We now want to discuss the general case. We will make use of the following claim: 
  
  \medskip
  
  {\bf Proposition 1.2.} Let $\Gamma $ be non-solvable, $F(\Gamma ) < \infty $, $\epsilon > 0$, and $x_0\in (0,1)$. Then there exist $\gamma _1, \gamma _2 \in [\Gamma ,\Gamma ]$ such that $0 < \gamma _1(x_0) - x_0 < \epsilon $ and $0 < x_0 - \gamma _2(x_0) < \epsilon $.

  \medskip
  
  Using this proposition, we immediately obtain the following lemma
  
  \medskip
  
  {\bf Lemma 1.3.} Let $\Gamma $ be non-solvable, $F(\Gamma ) < \infty $, $0 < a < b < 1$, $x_0\in (0,1)$. Then there exists $\gamma \in [\Gamma ,\Gamma ]$ such that $\gamma (x_0)\in (a,b)$. \ $\square $

  \medskip
  
  Proposition 1.2 follows from the following more general fact
  
  \medskip
  
  {\bf Proposition 1.4.} Let $\Gamma $ be non-solvable, $F(\Gamma ) < \infty $, $\epsilon > 0$, and $x_1, \ldots , x_n\in (0,1)$. Then there exists $\gamma \in [\Gamma ,\Gamma ]$ s.t.  $|\gamma (x) - x|< \epsilon , \forall x\in [0,1]$ and $Fix(\gamma )\cap \{x_1, \ldots , x_n\} = \emptyset $.
  
  \medskip
  
  {\bf Proof.} Let $F(\Gamma ) = N$. Let also $\eta \in \Gamma $ such that $\eta (x) < x$ near 1, and $max Fix(\eta ) < min \{x_1, \ldots , x_n\}$ (such $\eta $ exists by irreducibility of $\Gamma $). We will consider a finite family $\{\eta ^k \ : \ -nN-1 \leq k\leq nN+1\}$ of diffeomorphisms. By uniform continuity, there exists $\delta  > 0$ such that for all $x, y\in [0,1], -nN-1\leq k\leq nN+1$, we have  $|x-y| < \delta \Rightarrow |\eta ^k(x)-\eta ^k(y)| < \epsilon $.
  
  \medskip
  
  By Theorem B in [A], there exists $f\in [\Gamma ,\Gamma ]\backslash \{1\}$ such that $d(f) < \delta $. (Theorem B, as stated in [A], claims the existence of such $f$ in $\Gamma $, but it is immediately clear from the proof that $f$ can be chosen from the subgroup $[\Gamma , \Gamma ]$ \footnote {The proof of Theorem B relies on the proof of Theorem A. Notice that, in the proof of Theorem A, the maps $f$ and $g$ belong to the commutator subgroup. Hence the map $h_1^{-1}h_2$ at the end of the proof belongs to the commutator subgroup as well.}). Then $d(\eta ^{-k}f\eta ^k) < \epsilon $ for all $k\in \{-nN-1, -nN, \ldots , n, nN+1\}$. Moreover, by pigeonhole principle, there exists $k\in \{1, \ldots , nN, nN+1\}$ such that  $\{x_1, \ldots , x_n\}\cap Fix (\eta ^{-k}f\eta ^k) = \emptyset $. \ $\square $
  
   \bigskip
   
   We would like to emphasize that we do not know how to prove Proposition 1.2 directly, i.e. without using Proposition 1.4 which is significantly stronger. 
   
   \medskip
   
   Now we will prove a somewhat stronger claim which is more suitable for our purposes 
   
   \medskip
   
   {\bf Proposition 1.5.} Let $\Gamma $ be non-solvable, $F(\Gamma ) < \infty $, $\epsilon > 0$, and $0 < a < b < 1$. Then there exists $\gamma \in [\Gamma ,\Gamma ]$ s.t.  $|\gamma (x) - x|< \epsilon , \forall x\in [0,1]$ and $[a,b]\cap Fix(\gamma ) = \emptyset $.
   
   \medskip
   
   {\bf Proof.} In the proof of Proposition 1.4, it is sufficient to choose $\eta \in \Gamma $ such that $\eta ^{-1}(a) > b$. $\square $
   
   \medskip
   
   {\bf Remark 1.6.} Let us notice that, in Propositions 1.4 and 1.5, Theorem B (from [A]) is used only to guarantee the existence of $C_0$-small elements. But if $\Gamma \leq $Diff $ _{+}^2(I)$, then one can use Theorem C from [A] to guarantee the existence of $C_0$-small elements under a weaker condition that $\Gamma $ is not metaabelian. Hence, if $\Gamma \leq $Diff $ _{+}^2(I)$ then in Proposition 1.2, 1.4, 1.5 and Lemma 1.3 one can replace the condition "$\Gamma $ is not solvable" with "$\Gamma $ is not metaabelian". 
   
   \vspace{1cm}
   
   \begin{center} {\bf 2. Solvability of groups with small $N$} \end{center}
   
   \bigskip
   
   Our proof of the main theorem will use the assumption that $N\geq 5$. Therefore we need to take care of the cases $N\leq 4$ separately. It is also interesting to see how easy it is to prove Theorem 0.1 in the cases of $N\leq 4$, using the results of Section 1. 
   
   \medskip
   
  The cases of $N = 0, 1$ have been discussed earlier so it remains to study the cases of $N = 2, 3, 4.$ Our proofs for these cases are somewhat different from each other, and the techniques used in these proofs will be useful for the reader as a preparation for the proof of the main theorem.
   
   \medskip
   
   If $\gamma \in \Gamma , Fix(\gamma ) = \{a_1, \ldots , a_n\}$ where $a_1 < a_2 < \ldots < a_n$ then we will write $p_i(\gamma ) = a_i$ for all $1\leq i\leq n$. (so $p_i(\gamma )$ denotes the $i$-th fixed point of $\gamma $.)
   
   \medskip
   
   {\bf Proposition 2.1.} Let $\Gamma \leq G$ such that $F(\Gamma ) \leq 4$. Then $\Gamma $ is solvable.  
   
   \medskip
   
   {\bf Proof.} In $G$, the derived length of any solvable subgroup is universally bounded (see [N2]) therefore we may assume that $\Gamma $ is finitely generated with a fixed symmetric generating set. 
   
   \medskip
   
   Let $F(\Gamma ) = N$, and $N\in \{2,3,4\}$. Let also $f$ be the biggest generator of $\Gamma $. 
  
 \medskip
 
 Then $f\notin \Gamma _{\gamma }$ for all $\gamma \in \Gamma $. (in other words, $f$ is not infinitesimal w.r.t. to any element of $\Gamma $), and $[\Gamma , \Gamma ]\leq \Gamma _{f}$.
  
 \medskip
   
   Without loss of generality, we may assume that $f$ has no tangential fixed point. [Indeed, if $f$ has a tangential fixed point then if $\omega \in [\Gamma ,\Gamma ]$ is sufficiently close to the identity diffeomorphism in $C_0$-metric, and $Fix(\omega )\cap Fix(f) = \emptyset $ (such $\omega $ exists by Proposition 1.4), then there exists $\phi \in \{\omega f, \omega ^{-1}f\}$ such that $|Fix(\phi )|\geq |Fix(f)|$, moreover, either $\phi $ has no tangential fixed point or it has more fixed points than $f$. Notice also that $\phi $ is positive and $\phi \in \Gamma \backslash \Gamma _f$. Thus we can replace $f$ by $\phi $; and if we still have tangential fixed points we perform the same operation at most four times (since $F(\Gamma ) \leq 4$) until we find a positive element in $\Gamma \backslash \Gamma _f$ (which we still denote by $f$) without tangential fixed points.]

   \medskip
   
    Notice that $f(x) > x$ near zero. We may also assume that $f$ has at least two fixed points. (From the proof of Proposition 1.1, we see that, for a suitable choice of $h\in \Gamma _{f}$, and for sufficiently big $n$, the diffeomorphism $h^{-n}f_1$ has at least two fixed points where $f_1 = \omega f\omega ^{-1}$, for some $\omega \in \Gamma _f$. Hence, if necessary, we may replace $f$ with $h^{-n}f_1$). 
    
    \medskip
    
    a) Let $N=2$. Let also  $h\in \Gamma _{f}\backslash \{1\}$ such that $h$ has two fixed points (if such $h$ does not exist then $\Gamma _{f}$ is metaabelian by Solodov's Theorem, therefore $\Gamma $ is solvable with derived length at most 3). We may also assume that $h$ has no tangential fixed point, $Fix(h)\cap Fix(f) = \emptyset $, and $h(x) > x$ near zero.
 
 \medskip
    
   Let $Fix(h) = \{a, b\}, a < b$. Let us remind that $f$ has at least two fixed points, hence (since $N=2$) exactly two fixed points. 
   
   \medskip
 
  Let $c$ be the smallest fixed point of $f$. If necessary, using Lemma 1.3, by conjugating $f$ by an element of $\Gamma _f$, we may assume that $c>b$. Then, for sufficiently big $n$, the graphs of $f$ and $h^{n}$ will have at least two crossings on the interval $(0,a)$, and at least one crossing on the interval $(b,c)$. Hence $h^{-n}f$ will have at least three fixed points; contradiction.
   
    \medskip
    
   b) $N=3$. Let $h\in \Gamma _{f}\backslash \{1\}$ such that $h$ has three fixed points. (if such $h$ does not exist then then by part a), $\Gamma _{f}$ is solvable, therefore $\Gamma $ is solvable). We may also assume that $h$ has no tangential fixed points, $Fix(h)\cap Fix(f) = \emptyset $, and $h(x) > x$ near zero.
     
   \medskip
   
   Let $Fix(h) = \{a,b,c\}, a<b<c$.

   \medskip
   
   By assumption, $f$ has at least two fixed points. By conjugating $f$ if necessary, and using Lemma 1.3, we may assume that the smallest fixed point of $f$ (let us denote it by $p$) lies in $(b,c)$. Then, for sufficiently big $n\in \mathbb{N}$, the graphs of $f$ and $h^n$ have at least two crossings on $(0,a)$, and at least one crossing on $(b,p)$. If $f$ has three fixed points, then there will be one more crossings on $(p,1)$. If $f$ has two fixed points then there exists one more crossing either on $(p,c)$ (if the biggest fixed point of $f$ is less than $c$), or on $(c,1)$ (if the biggest fixed point of $f$ is bigger than $c$).  Contradiction. 
   
   \medskip
   
   c) Let now $N=4$. Again, there exists $h\in \Gamma _{f}\backslash \{1\}$ such that $h$ has four fixed points. (if such $h$ does not exist then by part b), $\Gamma _{f}$ is solvable, therefore $\Gamma $ is solvable). We may furthermore  assume that $h$ has no tangential fixed point, $Fix(h)\cap Fix(f) = \emptyset $, and $h(x) > x$ near zero.
     
   \medskip
   
   By assumption, $f$ has at least two fixed points. We may conjugate $h$ to $h_1$ by elements of $\Gamma _{f}$ such that $Fix(h_1)\subset (0, p_1(f))$. Then, for sufficiently big $n$,  $h_1^n$ and $f$ have at least three crossings on $(0, p_1(f))$. Hence $f_1 = h_1^{-n}f_1$ has at least three fixed points. (notice that $f_1$ is positive and $f_1\notin \Gamma _{f}.)$ 
 
  \medskip
  
  \begin{figure}[h!]
  \includegraphics[width=5in,height=5in]{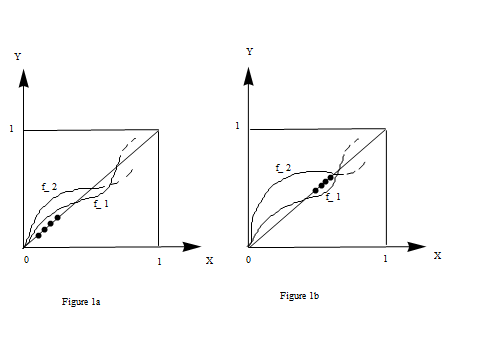}
\caption{the four dots in Figure 1a (Figure 1b) represent the fixed points of $h_2$ ($h_3$).}
\label{labelname}
\end{figure}
  
  \medskip
  
 We can conjugate $f$ to $f_2$ by elements of $\Gamma _{f}$ such that $p_1(f_2)\in (p_1(f_1), p_2(f_1))$. Then we can conjugate $h$ to $h_2$ by elements of $\Gamma _{f}$ such that $Fix(h_2)\subset (0, p_1(f_1))$ (Figure 1a). Then, conjugating $h_2$ by powers of $f_2$ we obtain $h_3$ where $Fix(h_3)\subset (p_1(f_1), p_1(f_2))$ (Figure 1b). Using Lemma 1.5, we can conjugate $f_1$ to $f_3$ by elements of $\Gamma _{f}$ such that the first two fixed points of $h_3$ lie in $(0, p_1(f_3))$ while the last two fixed points lie in $(p_1(f_3), p_2(f_3))$(Figure 2).     
   
  \medskip
  
   \begin{figure}[h!]
  \includegraphics[width=4in,height=4in]{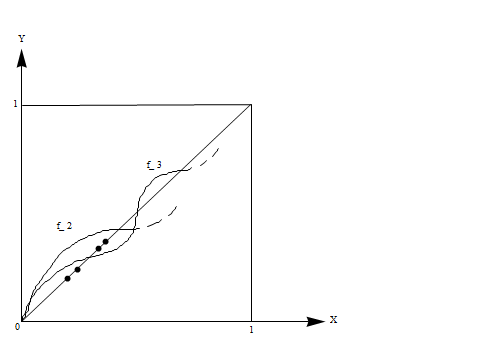}
\caption{the four dots represent the fixed points of $h_3$.}
\label{labelname}
\end{figure}
   
   \medskip
   
   Let $g_n = f_3^{-n}h_3f_3^n, n\in \mathbb{N}$. Notice that as $n\to \infty $, we have $p_1(g_n)\to 0, p_2(g_n)\to 0$ while  $p_3(g_n)\to p_2(f_3),  p_4(g_n)\to p_2(f_3)$. Then for sufficiently big $n$, using Proposition 1.3, $g_n$ can be conjugated to $\phi _n$ by an element of $\Gamma _{f}$ such that $\phi _n$ has two fixed points in $(0,p_1(f_3))$ and another two fixed points in $(p_2(f_3), p_3(f_3))$. Then, for sufficiently big $k$, the diffeomorphisms $\phi _n^{-k}$ and $f_3$ have at least three crossings on $(0,p_1(f_3))$, and at least two crossings on $(p_2(f_3), p_3(f_3))$. Contradiction.  
$\square $

\medskip

 {\bf Remark 2.2.} Interestingly, already for $N=5$, the elementary arguments in the proof of Proposition 2.1 do not seem to work, and one needs significantly new and more general ideas. (In fact, the value $N=5$ seems to be too big for ``elementary methods" and too small for ``general methods". We make extra efforts to cover this particular case). The essence of the problem is that Lemma 1.3. can be viewed as a dynamical transitivity of the action. However, we do not have a result about dynamical $k$-transitivity even when $k=2$. But as $N$ gets bigger, we need higher transitivity results to control the picture as we did in the proof of Proposition 2.1.  [{\em We call the action of $\Gamma $ on $I$ $k$-transitive if for all $z_1, \ldots , z_k \in (0,1)$, where $z_1 < z_2 < \ldots < z_k$ and for all open non-empty intervals $I_1, I_2, \ldots , I_k$ where $x < y$ if $x\in I_p, y\in I_q, 1\leq p < q\leq k$, there exists $\gamma \in \Gamma $ such that $\gamma (z_i)\in I_i, 1\leq i\leq k$. Proposition 1.5. is indeed significantly stronger than 1-transitivity.}]

 \vspace{2cm}
 
 \begin{center} {\bf 3. Proof of the Main Theorem} \end{center}
 
 \bigskip

   The following lemma is very interesting in itself. It will be used in the proof of Theorem 0.2, and the idea of the proof will be used repeatedly in the proof of Theorem 0.1.
   
   \medskip
   
   {\bf Lemma 3.1} Let $\Gamma $ be non-solvable, $F(\Gamma ) < \infty , f\in \Gamma , |Fix(f)| = n$. Then there exists $\psi \in [\Gamma ,\Gamma ]$ such that $|Fix(\psi )| \geq n-1$.  Moreover, $\psi $ can be taken arbitrarily close (in $C_0$-metric) to the identity.
   
   \medskip
   
   {\bf Proof.} Let $Fix(f) = \{c_1, \ldots , c_n\}$ where $c_1 < \ldots < c_n, \epsilon < \frac{1}{2}s(f), a = c_1 - \epsilon , b = c_n + \epsilon $ . By uniform continuity of $f$ and $f^{-1}$, there exists $\delta > 0$ such that if $d(h) < \delta $ then $d([f,h]) < \epsilon $. By Proposition 1.5, there exists $h \in [\Gamma , \Gamma ]$ such that $h(x) > x$ on $[a,b]$, $Fix(h) \cap [a,b] = \emptyset $, and $d(h) < \mathrm{min}\{\delta , \epsilon \}$.    
   
   \medskip
   
  Let now $g = h^{-1}fh$. Then $g$ has $n$ fixed points $d_1, \ldots , d_n$ where $d_1 < c_1 < d_2 < c_2 < \ldots < d_n < c_n$.
Then the element $\psi = f^{-1}g$ has at least $n-1$ fixed points $e_1, \ldots , e_{n-1}$ where $e_i\in (c_{i}, d_{i+1}), 1\leq i\leq n-1$. Moreover, since $d(h) < \delta $, we also have $d(\psi ) < \epsilon $. $\square $

  \medskip
  
  Now we want to observe the following simple lemma
  
  \medskip
  
  {\bf Lemma 3.2.} Let $f, h\in \mathrm{Homeo} _{+}(I), Fix(h)\cap [a,b] = \{a, b\}, Fix(f)\cap [a,b] = \emptyset $. Let also $f(a) \in (a,b)$ or $f(b)\in (a,b)$. Then the graphs of $h^{-1}fh$ and $f$ have at least one crossing on $(a,b)$. $\square $ 
  
  \medskip
  
  {\bf Remark 3.3.} Under the assumptions of the lemma (indeed, under even a slightly more general assumptions), $f$ and $h$ are called crossed elements. If $f, h$ are crossed then the subgroup generated by $f, h$ contains a non-abelian free semigroup([N1]). (we do not use this result here).
  
  \medskip
  
  In the proof of the main theorem, the case of $N=5$ requires more effort; we will need the following lemma in treating this case.
  
  \medskip 
  
  {\bf Lemma 3.4.} Let $\Gamma \leq G$ be a finitely generated non-solvable subgroup with a fixed symmetric generating set, $F(\Gamma ) = 5$, and $f$ be the biggest generator. Then at least one of the following conditions hold:
  
  \medskip
  
  (i) there exists $g\in \Gamma \backslash \Gamma _{f}$ with 5 fixed points.
  
  \medskip
  
  (ii) $\Gamma \backslash \Gamma _{f}$ contains a diffeomorphism with 4 fixed points and $\Gamma _f$ does not contain a diffeomorphism with 4 fixed points. \footnote {let us emphasize and clarify that when we say a diffeomorphism has $n$ fixed points, we mean that it has {\em exactly} $n$ fixed points in $(0,1)$; we often drop the phrase {\em exactly}}
  
  \medskip
  
  {\bf Proof.} Since $\Gamma $ is irreducible and $\Gamma _{f}$ is a normal subgroup, we may assume that $\Gamma _{f}$ has no global fixed point in $(0,1)$. As we already discussed, we may also assume that $f$ has at least two non-tangential fixed points. Assume that there is no $g\in \Gamma \backslash \Gamma _{f}$ with 5 fixed points. Then there exists $h\in \Gamma _{f}\backslash \{1\}$ with exactly five non-tangential fixed points and $h(x)>x$ near zero. Using Lemma 1.3, we can conjugate $h$ to $h_1$ by an element of $\Gamma _{f}$ such that the biggest fixed point of $h_1$ lies in $(p_1(f), p_2(f))$. We can also conjugate $h$ to $h_2$ such that $Fix(h_2)\subset (p_5(h_1), 1)$. [no need to use Lemma 1.3 here; we just use the fact $\Gamma _f$ has no global fixed point in $(0,1)$]. Then, we can conjugate $h_2$ to $h_3$ by powers of $h_1$ such that $Fix(h_3)\subset (p_5(h_1), p_2(f))$. Finally, using Lemma 1.3, we can conjugate $f$ to $f_1$ by an element of $\Gamma _{f}$ such that the first two fixed points of $h_3$ lie in $(0,p_1(f_1))$ and the last three fixed points lie in $(p_1(f_1), p_2(f_1))$. Then, for sufficiently big $n$,  the diffeomorphisms $f_1$ and $h_3^n$ have at least three crossings in $(0,p_1(f_1))$, and at least one crossing in $(p_1(f_1), p_2(f_1))$. Thus, for sufficiently big $n$, the diffeomorphism $f_2 = h_3^{-n}f_1$ has at least four fixed points, hence exactly four fixed points. Notice that $f_2\notin \Gamma _{f}$ and $f_2$ is positive.
  
  \medskip
  
 Now, if there is no diffeomorphism in $\Gamma _f$ with four fixed points then we are done. But if such a diffeomorphism exists then we fall in the case of the proof of Proposition 2.1, part c), (i.e. having a diffeomorphism $f_2\in \Gamma \backslash \Gamma _f$ with at least three fixed points and a diffeomorphism $\Gamma _f$ with four fixed points) thus obtain a diffeomorphism with five fixed points in $\Gamma \backslash \Gamma _f$. $\square $
    
\bigskip
  
  Now we are ready to prove the main result.
  
  \medskip
  
  {\bf Proof of Theorem 0.1.} We will assume that $N\geq 5$. Let $g\in \Gamma $ with $N$ fixed points, and $g(x) > x$ near zero. We may assume that none of the fixed points of $g$ is tangential. By conjugating $g$ we obtain $g_1$ and $g_2$ such that $Fix(g_1)\subset (0, \frac{1}{10}), Fix(g_2)\subset (\frac{9}{10},1)$. Let $Fix(g_1) = \{a_1, \ldots , a_N\}, Fix(g_2) = \{b_1, \ldots , b_N\}$ where the elements are listed from the least to the biggest.
  
  \medskip
   
  Using Proposition 1.5, as in the proof of Lemma 3.1, we obtain $h_1\in \Gamma $ such that $g_2$ and $h_1^{-1}g_2h_1$ have at least one crossing in each of the intervals $(b_i, b_{i+1}), 1\leq i\leq N-1$, moreover, if $d(h_1)$ is sufficiently small, then $d(f) < \mathrm{min}\{\frac{1}{10},\frac{1}{2}s(g_1)\}$ where $f = g_2^{-1}h_1^{-1}g_2h_1$. Thus the diffeomorphism $f$ has at least $N-1$ fixed points in $(b_1, 1)\subset (\frac{9}{10},1)$. Then $f$ has at most one fixed point in $(0, \frac{9}{10})$.
  
  \medskip
  
  Using Proposition 1.5 again, we obtain $h_2$ such that the diffeomorphisms $g_1$ and $h_2^{-1}g_1h_2$ have at least one crossing in each of the intervals $(a_i, a_{i+1}), 1\leq i\leq N-1$, moreover, if $d(h_2)$ is sufficiently small, then $d(h) < \mathrm{min}\{\frac{1}{10},\frac{1}{2}s(f)\}$ where $h = g_1h_2^{-1}g_1h_2$. Thus the diffeomorphism $h$ has at least $N-1$ fixed points in $(0, \frac{1}{10})$.
  
  \medskip
  
  Let us now assume that $N\geq 9$. Then (let us recall that $d(f) < \mathrm{min}\{\frac{1}{10},\frac{1}{2}s(g_1)\}$), if $d(h_2)$ is sufficiently small, there exists at least four mutually disjoint non-empty intervals $I_i = (p_i,q_i), 1\leq i\leq 4$ in $(0,\frac{1}{5})$ such that for all $i\in \{1,2,3,4\}$ 
  
  \medskip
  
  (a1) $p_i, q_i\in Fix(h)$;
  
  \medskip
  
  (a2) $Fix(h) \cap I_i = \emptyset $;
  
  \medskip
  
  (a3) $f(p_i) \in I_i $ or $f(q_i) \in I_i$. 
  
  \medskip
  
  Since $f$ has at most one fixed point in $(0,\frac{9}{10})$, $f$ has no fixed points in at least three of these intervals. Thus, applying Lemma 3.2, we obtain that $f$ and $h^{-1}fh$ have at least three crossings in $(0,\frac{9}{10})$. But since $d(h) < \frac{1}{2}s(f)$, we have at least one crossing in each of the intervals $(c_j, c_{j+1}), 1\leq j\leq M-1$ where $c_1, \ldots , c_{M}$ are all fixed points of $f$ listed from least to the biggest, and $M\in \{N-1,N\}$ (we know that $f$ has at least $N-1$ fixed points). Hence, $f$ and $h^{-1}fh$ have at least $3+(N-2) = N+1$ crossings. Contradiction.   
   
  \medskip
  
  When $N\in \{5,6,7,8\}$ we need to sharpen our observation. Indeed, for $N\in \{6,7,8\}$, we have either 
  
  {\em Case 1.} there exist four mutually disjoint non-empty intervals $I_i = (p_i,q_i), 1\leq i\leq 4$ in $(0,\frac{1}{5})$ such that for all $i\in \{1,2,3,4\}$ the conditions (a1)-(a3) hold, 
  
  or 
  
  {\em Case 2.} there exist three mutually disjoint non-empty intervals $I_i = (p_i,q_i), 1\leq i\leq 3$ in $(0,\frac{1}{5})$ such that for all $i\in \{1,2,3\}$ the conditions (a1)-(a3) hold, moreover, in two consecutive intervals, either $h(x) > x$ on both of these intervals, or $h(x) < x$.
  
  \medskip
  
  In both cases, we obtain that $f$ and $h^{-1}fh$ have at least three crossings in $(0,\frac{9}{10})$, hence, at least $N+1$ crossings total; contradiction.
  
  \medskip
  
  The case of $N=5$ needs a special treatment. We may assume that $\Gamma $ is finitely generated with a fixed generating set and let $\theta $ be the biggest generator. 
  
  \medskip
  
  First, let us assume that $\Gamma \backslash \Gamma _{\theta }$ contains a diffeomorphism with 5 fixed points. The we choose $g$ to be this diffeomorphism and we let $g_1, g_2$ be obtained from $g$ by conjugating by elements from $\Gamma _{\theta }$.  
  
  \medskip
  
   Then $f\in \Gamma _{\theta }$. Hence if $f$ has a fixed point in $(0, b_1)$ then for sufficiently big $n$, $f^n$ and $g_2$ will have at least two crossings in $(0,b_1)$, hence, at least $(N-1)+2$ crossings total, contradiction. Thus, we may assume that $f$ has no fixed point in $(0,b_1)\supset (0, \frac{9}{10})$, and the conclusion follows as in the case of $N\geq 6$. 
   
   \medskip
   
   If $\Gamma \backslash \Gamma _{\theta }$ does not contain a diffeomorphism with 5 fixed points, then let $\Gamma _1$ be an arbitrary finitely generated non-solvable subgroup of $\Gamma _{\theta }$ with a fixed finite symmetric generating set, and let $\eta \in \Gamma _1$ be the biggest generator. By Lemma 3.4, the normal subgroup $(\Gamma _1)_{\eta}$ contains no diffeomorphism with four fixed points. Hence, since $\Gamma _1$ does not contain a diffeomorphism with four fixed points, again by Lemma 3.4, we obtain that $\Gamma _1\backslash (\Gamma _1)_{\eta }$ contains a diffeomorphism with five fixed points. Thus applying the previous argument to this case we obtain a contradiction. Hence $\Gamma _1$ is solvable. But since $\Gamma _1$ is an arbitrary finitely generated subgroup of $\Gamma _{\theta }$, we obtain that $\Gamma _{\theta }$ is solvable. Hence $\Gamma $ is solvable.  $\square $
 
 \medskip
 
 For the purposes of clarity, let us emphasize that the proof of Theorem 0.1. does not use an inductive argument on $N\geq 5$.

 \bigskip
 
 Now we are ready to prove Theorem 0.2. Taking into account Remark 1.6, the proofs of Proposition 2.1. (for $N\in \{2,3,4\}$) and of Theorem 0.1. (for $N\geq 6$) go {\em mutatis mutandis} to obtain a contradiction under the assumptions that $\Gamma \leq $ Diff $_{+}^{2}(I)$ is non-metaabelian and $F(\Gamma ) = N$. \footnote{We also need to show that the claims in parenthetical remarks in parts a)-c) of the proof of Proposition 2.1 hold as well. More precisely, we need to show that if $N\in \{2,3,4\}$ then $\Gamma _f$ contains a diffeomorphism with $N$ fixed points. Indeed, for $N = 2$, if $\Gamma _f$ does not have such a diffeomorphism then by Theorem 1.5 of [FF], there exists an element of $h\in \mathrm{Homeo}_{+}(I)$ such that $h\Gamma _fh^{-1}$ is a subgroup of $\mathrm{Aff}_{+}I$ (this latter group, by definition, is obtained from $\mathrm{Aff}_{+}\mathbb{R}$ by conjugation by a tangent function). Then the subgroup $\Gamma _1 = \{g\in h\Gamma _fh^{-1} \ | \ g = id \ \mathrm{or} \ Fix(g) = \emptyset \}$ contains elements arbitrarily close to the identity in $C_0$-metric. Let $\phi \in h\Gamma h^{-1}$ with two fixed points. Then there exists $\omega \in \Gamma _1$ such that $[\omega , \phi ] = \omega \phi \omega ^{-1}\phi ^{-1}$ has one fixed point. But $\phi \omega ^{-1}\phi ^{-1}$ has no fixed point, hence it belongs to $\Gamma _1$. Then $[\omega , \phi ]\in\Gamma _1$, hence it has no fixed point; contradiction. For the case of $N = 3$, assume $\Gamma _f$ has no element with three fixed points. If it also has no element with two fixed points then we run the same argument. But if $F(\Gamma _f) = 2$ then by part a), we obtain a contradiction.  Similarly, for the case of $N = 4$: assume $\Gamma _f$ has no element with 4 fixed points. If $F(\Gamma _f) = 1$ then we run the previous argument; but if $F(\Gamma _f)\in \{2,3\}$ then by parts a) and b) we obtain a contradiction again.} 
 
 The case of $N=5$ is again slightly subtle. In this case, we can only conclude that the group $\Gamma _1$ (at the end of the proof of Theorem 0.1.) is metaabelian. Since $\Gamma _1$ is an arbitrary finitely generated subgroup of $\Gamma _{\theta }$ then we can conclude that $\Gamma _{\theta }$ is metaabelian. This does not in automatically imply that $\Gamma $ is metaabelian (it only implies that $\Gamma $ is solvable of derived length at most 3). So let us assume that $\Gamma $ is not metaabelian. 
 
 \medskip
 
 If $\Gamma \backslash \Gamma _{\theta }$ contains a diffeomorphism with 5 fixed points then the claim again follows from the proof of Theorem 0.1;  we choose $g$ to be a diffeomorphism from $\Gamma \backslash \Gamma _{\theta }$ with five fixed points, and we let $g_1, g_2$ be obtained from $g$ by conjugating by elements from $\Gamma _{\theta }$.
 
 \medskip
 
 If $\Gamma \backslash \Gamma _{\theta }$ does not contain a diffeomorphism with five fixed points then there exists $\xi \in \Gamma _{\theta }$ with five fixed points. 
  
 \medskip
 
 On the other hand, since $\Gamma $ is not metaabelian, $[\Gamma , \Gamma ]$ contains a non-trivial element $\eta $ with $d(\eta ) < \frac{1}{2}s(\xi ) $ and $Fix(\eta ) \cap Fix(\xi ) = \emptyset $. Then $Fix(\eta \xi \eta ^{-1})\cap Fix(\xi ) = \emptyset $ and    $\gamma = \eta \xi \eta ^{-1}\xi ^{-1}$ has at least four fixed points where $\gamma \neq 1$.

 \medskip
 
  Since $[\Gamma , \Gamma ]\leq \Gamma _{\theta }$ we have $\eta \in \Gamma _{\theta}$. But we also have $\xi \in \Gamma _{\theta }$ therefore $\gamma \in [\Gamma _{\theta }, \Gamma _{\theta }]$. 
  
  \medskip
  
  Now, since $\Gamma $ is non-metaabelian, there exists $\omega \in [\Gamma , \Gamma ]\leq \Gamma _{\theta }$ such that $d(\omega ) < \frac{1}{2}s(\gamma )$ and $Fix(\omega )\cap Fix(\gamma ) = \emptyset $. Then $Fix(\omega \gamma \omega ^{-1})\cap Fix(\gamma ) = \emptyset \ (\star )$. Thus we obtain that $\omega \gamma \omega ^{-1}$ and $\gamma $ do not commute (if they do then by $(\star )$, this forces $\gamma $ to have infinitely many fixed points). But $\gamma \in [\Gamma _{\theta }, \Gamma _{\theta }]$, therefore $\omega \gamma \omega ^{-1}\in [\Gamma _{\theta }, \Gamma _{\theta }]$. On the other hand, $[\Gamma _{\theta }, \Gamma _{\theta }]$ is Abelian. Then $[\gamma , \omega \gamma \omega ^{-1}] = 1$. Contradiction. $\square $

  \vspace{1cm}
  
  {\bf References}
  
  \bigskip
  
  [A] Akhmedov A. \ A weak Zassenhaus lemma for subgroups of Diff(I). To appear in {\em Algebraic and Geometric Topology}. \\ http://arxiv.org/pdf/1211.1086.pdf
  
  \medskip
  
  [B] T.Barbot, \ Characterization des flots d'Anosov en dimension 3 par leurs feuilletages faibles, \ {\em Ergodic Theory and Dynamical Systems} {\bf 15} (1995), no.2, 247-270.
  
  \medskip
  
  [FF] B.Farb, J.Franks, \ Groups of homeomorphisms of one-manifolds II: Extension of H\"older's Theorem. \ {\em Trans. Amer. Math. Soc.} {\bf 355} (2003) no.11, 4385-4396.
  
  \medskip
  
  [K] N.Kovacevic, \ M\"obius-like groups of homeomorphisms of the circle. \ {\em Trans. Amer. Math. Soc.} {\bf 351} (1999), no.12, 4791-4822.
  
  \medskip
  
  [N1] Navas, A. Groups of Circle Diffemorphisms. \ Chicago Lectures in Mathematics, 2011. {\em http://arxiv.org/pdf/0607481} 
  
  \medskip
  
  [N2] Navas, A. Growth of groups and diffeomorphisms of the interval. \ {\em Geom. Funct. Anal.} \ {\bf 18} 2008, {\bf no.3}, 988-1028.

 \end{document}